\documentclass[11pt]{amsart}
\usepackage{amssymb}
\usepackage{amsfonts}
\usepackage{color}

\usepackage{diagrams}

\newarrow{Dash}{}{dash}{}{dash}>
\newarrow{Into}C--->
\newarrow{Onto}----{>>}
\newarrow{Equal}===={=}
\newarrow{Dots}....>

\newtheorem{theorem}{Theorem}[section]
\newtheorem{proposition}[theorem]{Proposition}

\newtheorem{definition}[theorem]{Definition}

\newtheorem{example}[theorem]{Example}

\theoremstyle{remark}

\renewenvironment{proof}{{\noindent\bf Proof.}}{\hfill $\Box$\par\vskip3mm}

\newcommand{\Hom}{{\rm Hom}}

\newcommand{\Cc}{\mathcal{C}}
\newcommand{\Dd}{\mathcal{D}}

\newcommand{\Ff}{\mathcal{F}}

\newcommand{\Ll}{\mathcal{L}}
\newcommand{\Mm}{\mathcal{M}}
\newcommand{\Pp}{\mathcal{P}}

\newcommand{\Ss}{\mathcal{S}}

\def\NN{{\mathbb N}}
\def\CC{{\mathbb C}}

\begin{document}
\title{When is $\prod$ isomorphic to $\bigoplus$}

\begin{abstract}
For a category $\Cc$ we investigate the problem of when the coproduct $\bigoplus$ and the product functor $\prod$ from $\Cc^I$ to $\Cc$ are isomorphic for a fixed set $I$, or equivalently, when the two functors are Frobenius functors. We show that for an {\bf Ab} category $\Cc$ this happens if and only if the set $I$ is finite. Moreover, this happens even in a much more general case, if there is a morphism in $\Cc$ that is invertible with respect to the addition of morphisms. If $\Cc$ does not have this property then we give an example to see that the two functors can be isomorphic for infinite sets $I$. However we show that $\bigoplus$ and $\prod$ are always isomorphic on a suitable subcategory of $\Cc^I$ which is isomorphic to $\Cc^I$ but is not a full subcategory. For the module category case we provide a different proof to display an interesting connection to the notion of Frobenius corings.
\end{abstract}

\author{Miodrag Cristian Iovanov}
\thanks{$^*$ Partially supported by the Flemish-Romanian project "New Techniques in Hopf algebras and Graded Rings Theory" (2005) and by the BD type grant BD86 2003-2005 by CNCSIS.}
\date{}
\maketitle

\section*{Introduction}
\noindent

Let $\Cc$ be a category and denote by $\Delta$ the diagonal functor from $\Cc$ to $\Cc^I$ taking any object $X$ to the family $(X)_{i\in I}\in\Cc^I$. Recall that $\Cc$ is an $\mathbf {Ab}$ category if for any two objects $X,Y$ of $\Cc$ the set ${\rm Hom}(X,Y)$ is endowed with an abelian group structure that is compatible with the composition. We shall say that a category is an $\mathbf{AMon}$ category if the set ${\rm Hom}(X,Y)$ is (only) an abelian monoid for every objects $X,Y$. Following \cite{McL}, in an $\mathbf{Ab}$ category if the product of any two objects exists then the coproduct of any two objects exists and they are isomorphic. Moreover, if this is the case, the product $M\times N$ with projections $p:M\times N\rightarrow M$ and $q:M\times N\rightarrow N$ and injections $i:M\rightarrow M\times N$ and $j:N\rightarrow M\times N$ defined such that $pi=1_M$, $qi=0$, $pj=0$, $qj=1_N$ form the bi-product of $M$ and $N$. This is also true if we only assume that $\Cc$ is an $\mathbf{AMon}$ category. \\
Given a family of objects $(X_i)_{i\in I}$ that has product $P$ and coproduct $C$, we may ask the question of when are they isomorphic. As there are very large classes of examples in which this is true (as it will be shown), the problem should be put in a functorial manner. If $I$ is a set, following \cite{McL}, a coproduct $\coprod\limits_{i\in I}{X_i}$ exists for every family of objects $(X_i)_{i\in I}$ if and only if the diagonal functor $\Delta$ las a left adjoint which can be constructed by considering the association $(X_i)_{i\in I}\longrightarrow\coprod\limits_{i\in I}{X_i}$. Dually, a product will be a right adjoint to $\Delta$. It is then natural to put the problem in this context, that is, when the product and coproduct functors are isomorphic functors. This is equivalent to asking when $\Delta$ has the same right and left adjoint, that is, when is it a Frobenius functor (or equivalently, when the product and/or coproduct are Frobenius; see \cite{Mo}, \cite{CMZ}). \\
One may expect that if the coproduct of a family $(X_i)_{i\in I}$ exists, then coproducts of subfamilies exist, and then the coproduct of a void family (of objects), which is an initial object. Thus we will assume that a zero object exists in $\Cc$, that is, an initial and terminal object. Then the existence of coproducts indexed by a set $I$ implies the existence of coproducts indexed by sets of cardinality less or equal to $I$, by taking the coproduct of $(X_j)_{j\in J}$ to be the coproduct of $(X_j)_{j\in J}\cup (0)_{i\in I\setminus J}$ for a subset $J$ of $I$. It is proved in \cite{Mc1} that in a category with a null object, finite products and finite coproducts and such that the canonical morphism $M\coprod N\rightarrow M\times N$ is always an isomorphism, every set ${\rm Hom}(M,N)$ can be endowed with an addition $f+g=\varepsilon(f\times g)\eta$ (where $\eta:M\rightarrow M \times M$ and $\varepsilon:N\coprod N\rightarrow N$ are the canonical morphisms) which is compatible with the composition of morphisms and makes ${\rm Hom}(M,N)$ into a commutative monoid. Thus the minimal setting for the problem will be that of an $\mathbf{AMon}$ category $\Cc$ that has a zero object, as it will be shown that if the product and coproduct functors indexed by a set $I$ are isomorphic in a category $\Cc$ with a zero object, then the canonical morphism from the coproduct to the product is an isomorphism, and then, by using the zero object, it follows that the isomorphism also holds for finite families of objects, making $\Cc$ into an $\mathbf{AMon}$ category. It will be shown that if in an $\mathbf{AMon}$ category with a zero object there is a nonzero morphism $M\stackrel{f}{\rightarrow}N$ that is invertible in the abelian monoid ${\rm Hom}(M,N)$ then the product and coproduct functors indexed by a set $I$ are isomorphic if and only if the set $I$ is finite. If such a morphism does not exist, the two functors can be isomorphic for infinite sets $I$. \\
In the second part a different solution is presented for the module category case. Although some quite direct considerations lead straight to the result, we wish to evidence an interesting connection to the theory of corings and comodules. Though introduced initially by Sweedler, corings and then categories of comodules over corings have proved to be generalizations of various categories such as the ones of Hopf modules, Yetter Drinfel'd modules, Doi Koppinen modules and entwined modules, categories of modules and comodules, of graded modules and categories of chain complexes, providing a unifying context for all these structures. The category of modules can be regarded as a category of comodules over the Sweedler canonical coring. Straightforward general considerations allow the viewing of the product category $\Mm_R^I$ as a category of comodules over a certain $R$ coring $C$ and the coproduct functor as the forgetful functor from this category to that of $R$ modules. Then the fact that the coproduct is a Frobenius functor will equivalently translate to the fact that the coring $C$ is Frobenius, which by a finiteness theorem will imply that $C$ must be finitely generated and the index set must be finite.

We will make use of the following remark: 
Let $\Cc$ be an $\mathbf{AMon}$ category that has a zero object and finite canonically isomorphic products and coproducts. Let $M,N$ be two objects of $\Cc$. Denote by $M\times M$ the (bi)product of $M$ by $M$ with $p,q:M\times M\rightarrow M$ the canonical projections, $N\times N$ the coproduct (biproduct) of $N$ by $N$ with $i,j:N\rightarrow N\times N$ the canonical injections, $\eta:M\rightarrow M\times M$ the unique morphism such that $p\eta=q\eta=1_M$ (the unit of the adjunction of the diagonal functor $\Cc\rightarrow \Cc\times\Cc$ with the product functor) and $\varepsilon:N\times N\rightarrow N$ the unique morphism such that $\varepsilon i=\varepsilon j=1_N$ (the counit of the adjunction of the coproduct functor with the diagonal functor $\Cc\rightarrow \Cc\times\Cc$). Then for two morphisms $f,g:M\rightarrow N$ we have $f+g=\varepsilon(ifp)\eta+\varepsilon(jgq)\eta=\varepsilon(f\times 0)\eta+\varepsilon(0\times g)\eta=\varepsilon(f\times g)\eta$.


\section{The General Case}

Let $I$ be a fixed set and $\Cc$ be an $\mathbf{AMon}$ category with a zero object $0$ and such that for any family $(X_i)_{i\in I}$ of objects of $\Cc$ there is a coproduct $\bigoplus\limits_{I}X_i$ and a product $\prod\limits_{I}X_i$ with the usual universal properties. Equivalently, $\bigoplus\limits_{I}$ is a left adjoint to $\Delta$ and $\prod\limits_{I}$ is a right adjoint to the diagonal functor $\Delta$. For objects $A,B$ denote by $0_{A,B}\in\rm{Hom}(A,B)$ (or simply $0$) the composition $A\stackrel{0}{\rightarrow}0\stackrel{0}{\rightarrow}B$; bythen $0_{A,C}=0_{B,C}\circ 0_{A,B}$ for every objects $A,B,C$ in $\Cc$ and $f\circ 0=0$; $0\circ g=0$ for every morphisms $f,g$ in $\Cc$.
We denote by $\sigma$ the unit of the adjunction $(\bigoplus,\Delta)$ and by $p$ the counnit of the adjunction $(\Delta,\prod)$. Let $A=(A_i)_{i\in I}$ be an object of $\Cc^I$. We will write shortly $\sigma_k$ (respectively $p_k$) for $\sigma_A^k$ (respectively $p_A^k$) - the $k$-th component of $\sigma_A$ (respectively $p_A$). For a functorial morphism $\Theta:\bigoplus\limits_{i\in I}\rightarrow\prod\limits_{i\in I}$ we denote by $\Theta_A^{k,l}=p_l\Theta_A\sigma_k$. The canonical functorial morphism $F_A:\bigoplus\limits_{i\in I}A_i\rightarrow\prod\limits_{i\in I}A_i$ is (by definition) the one with $F_A^{k,l}=0$ when $k\neq l$ and $F_A^{k,k}=1_{A_k},\,\forall k$. Also denote by $\Gamma_k:\Cc\rightarrow \Cc^I$ the functor associating to any $X\in \Cc$ the family $(X_i)_{i\in I}$ having $X_k=X$ and $X_i=0$ for $i\neq k$. 

\begin{proposition}\label{gen1}
A functorial morphism $\Theta:\bigoplus\limits_{i\in I}\rightarrow\prod\limits_{i\in I}$ is uniquely determined if for every object $X\in \Cc$ a morphism $\Lambda_X=(\Lambda_{X}^{i})_{i\in I}\in {\rm Hom}_{\Cc^I}(\Delta(X),\Delta(X))$ is given such that:
\begin{enumerate}
\item[(i)] $\Theta_X^{k,l}=0$ for $k\neq l$;
\item[(ii)] $\Theta_X^{k,k}=\Lambda_{X_k}^k$;
\item[(iii)] $f\Lambda_X^k=\Lambda_Y^lf$ for any objects $X,Y$ of $\Cc$ and any morphisms $f\in {\rm Hom}(X,Y)$.
\end{enumerate}
Moreover, if $\Theta$ is a functorial isomorphism then $\Lambda_X^i$ are isomorphisms.
\end{proposition}
\begin{proof}
If $A=(A_i)_{i\in I}$ and $B=(B_i)_{i\in I}$ are families of objects of $\Cc$ (that is, objects of $\Cc^I$) then we have the following commutative diagram
\begin{diagram}
& A_{k} & \rTo^{\sigma_{k}} & \bigoplus\limits_{i}A_{i} & \rTo^{\Theta_{(A_{i})_{i}}} & \prod\limits_{i}A_i & \rTo^{p_l} & A_l \\
(*) & \dTo^{f_k}& & \dTo^{\bigoplus\limits_{i}f_i} & & \dTo^{\prod\limits_{i} f_i} & & \dTo_{f_l} \\ 
& B_{k} & \rTo^{\sigma_{k}} & \bigoplus\limits_{i}B_{i} & \rTo^{\Theta_{(B_{i})_{i}}} & \prod\limits_{i}B_i & \rTo^{p_l} & B_l
\end{diagram}
where the left and right diagrams commute because of the naturality of $\sigma$ and $p$ and the middle diagram commutes by the naturality of $\Theta$. Thus we obtain $f_l\circ\Theta_A^{k,l}=\Theta_B^{k,l}\circ f_k$ for any $f_k\in {\rm Hom}_{\Cc}(A_k,A_k)$ and $f_l\in {\rm Hom}_{\Cc}(A_l,A_l)$, which yields (i) by taking $f_l=1_{A_l}$ and $f_k=0$, and together with (ii), it implies (iii). Denote by $\Lambda_X^k=\Theta_{\Gamma_k(X)}^{k,k}$. Taking $B=\Gamma_k(A_k)$ and $f_k=1_{A_k}$, $f_i=0:A_i\rightarrow 0$ in the upper diagram we see that $\Theta_A^{k,k}=\Lambda_{A_k}^k$, which only depends on $A_k$. It is obvious, by the universality of $\bigoplus$ and $\prod$ that $\Theta_A$ is uniquely determined by the morphisms $\Theta_A^{k,l}$. For the last part, consider the commutative diagram 
\begin{diagram}
\bigoplus\Gamma_k(X) & \rTo^{\Theta_{\Gamma_k(X)}} & \prod\Gamma_k(X) \\
\uTo^{\sigma_{\Gamma_k(X)}} & & \dTo_{p_{\Gamma_k(X)}} \\
X & \rTo^{\Lambda_X^k} & X \\
\end{diagram}
where it is easy to see that the vertical arrows are isomorphisms.
\end{proof}


\begin{proposition}\label{p1}
If the functors $\bigoplus$ and $\prod$ are isomorphic then the canonical morphism from $\bigoplus$ to $\prod$ is a functorial isomorphism. 
\end{proposition}
\begin{proof}
Let $\Theta$ be an isomorphism between $\bigoplus$ and $\prod$ and $\Lambda$ the morphism from Proposition \ref{gen1}. For each $X\in \Cc$, denote $L_X^k=(\Lambda_X^k)^{-1}$ and for every object $A=(A_{i})_{i\in I}\in \Cc^I$ consider the morphism $L_{A}=\bigoplus\limits_{i\in I}L_{A_i}^i:\bigoplus\limits_{i\in I}A_i\rightarrow\bigoplus\limits_{i\in I}A_i$. We see that this is a functorial morphism, as we have $\bigoplus\limits_{i\in I}L_{A_i}^i\circ\bigoplus\limits_{i\in I}f_i=\bigoplus\limits_{i\in I}(L_{A_i}^i\circ f_i)=\bigoplus\limits_{i\in I}(f_i\circ L_{A_i}^i)=\bigoplus\limits_{i\in I}f_i\circ\bigoplus\limits_{i\in I}L_{A_i}^i$ by (iii) in Proposition \ref{gen1} and the functoriality of $\bigoplus$, so we have the commutative diagram:
\begin{diagram}
\bigoplus\limits_{i\in I}A_{i} & \rTo^{\bigoplus\limits_{i\in I}L_{A_i}^i} & \bigoplus\limits_{i\in I}A_i \\
\dTo^{\bigoplus\limits_{i\in I}f_i} & & \dTo_{\bigoplus\limits_{i\in I}f_{i}} \\
\bigoplus\limits_{i\in I}A_{i} & \rTo_{\bigoplus\limits_{i\in I}L_{A_i}^i} & \bigoplus\limits_{i\in I}A_i\\
\end{diagram}
Now taking $f_i=L_{A_i}^i$ in diagram (*) we obtain $p_l\Theta_AL_A\sigma_k=p_l\Theta_A\sigma_kL_{A_k}^k$ and therefore $p_l\Theta_AL_A\sigma_k=0$ if $k\neq l$ (by Proposition \ref{gen1} (i)) and gives $p_k\Theta_AL_A\sigma_k=\Lambda_{A_k}^k\circ L_{A_k}^k=1_{A_k}$ for $k=l$, showing that $\Theta\circ L$ is the canonical morphism from $\bigoplus$ to $\prod$ so it is an isomorphism, as $L$ is an isomorphism with inverse $L'_A=\bigoplus\limits_{i\in I}\Lambda_{A_i}^i$.
\end{proof}

If $\Cc$ is a category with null (zero) object and the functors $\bigoplus\limits_{i\in I}$ and $\prod\limits_{i\in I}$ are isomorphic, then, by Proposition \ref{p1} they must be isomorphic by the canonical functor $\bigoplus\limits_{i\in I}\rightarrow\prod\limits_{i\in I}$. If $I$ is infinite then we can find an injection of sets $\NN\subseteq I$. For an object $X$ of $\Cc$ consider the family $\overline{X}=(X_i)_{i\in I}$ with $X_i=X$ for $i\in \NN$ and $X_i=0$ for $i\in I\setminus \NN$. It is not difficult to see that the fact that the canonical morphism for this family is an isomorphism implies the fact that the canonical morphism from $X^{(\NN)}=\bigoplus\limits_{n\in \NN}X$ to $X^{\NN}=\prod\limits_{n\in \NN}X$ is an isomorphism. By a similar argument we can show that the canonical morphism from $X\oplus Y$ to $X\times Y$ is an isomorphism, thus for any two objects we have a biproduct of $X$ and $Y$, so $\Cc$ is an {\bf AMon}-category by the remarks in the Introduction.

\begin{proposition}\label{p2}
Let $\Cc$ be an {\bf AMon}-category with zero object. Suppose the set $I$ is infinite and the canonical functor from $\bigoplus\limits_{i\in I}$ to $\prod\limits_{i\in I}$ is an isomorphism. Then no nonzero morphism $f:X\rightarrow Y$ in $\Cc$ is invertible in the monoid ${\rm Hom}(X,Y)$. 
\end{proposition}

\begin{proof}
Let $A$ be an object of $\Cc$. Denote by $r_0,r_1:A^\NN\times A^\NN\rightarrow A^\NN$ the canonical projections of $A^\NN\times A^\NN$, by $s_0,s_1:A^{(\NN)}\rightarrow A^{(\NN)}\times A^{(\NN)}$ the canonical injections of $A^{(\NN)}\times A^{(\NN)}$ ($\times$ is the biproduct in $\Cc$), $p_n$ (respectively $\sigma_n$) the canonical projections of $A^{\NN}$ (respectively injections of $A^{(\NN)}$). \\
Now consider the following diagram in $\Cc$:
\begin{diagram}
A & & \rTo^\Pi & & A^{\NN} & \pile{\lTo^\Lambda\\ \rTo_\Gamma} & A^{(\NN)} & & \rTo^\Sigma & & A \\
\dEqual & & (1) & & \uTo^{v} & (2) & \dTo_{u} & & (3) & & \dEqual \\
A & \rTo_{\Pi} & A^{\NN} & \rTo_{\eta} & A^{\NN}\times A^{\NN} & \pile{\lTo^{\Lambda\times\Lambda} \\ \rTo_{\Gamma\times\Gamma}} & A^{(\NN)}\times A^{(\NN)} & \rTo_{\varepsilon} & A^{(\NN)} & \rTo_{\Sigma} & A \\
\end{diagram}
where the morphisms are defined uniquely by the following relations (all canonical): 
\begin{eqnarray}
\eta & : & r_0\eta=r_1\eta=1_A \label {e1}\\
\varepsilon & : & \varepsilon s_0=\varepsilon s_1=1_A \label{e2}\\
\Sigma & : & \Sigma \sigma_n=1_A \label{e3}\\
\Pi & : & p_n\Pi=1_A \label{e4}\\
u & : & u\sigma_{2n}=s_0\sigma_n \label{e5}\\
 & & u\sigma_{2n+1}=s_1\sigma_n \label{e6}\\
v & : & p_{2n}v=p_nr_0 \label{e7}\\
 & & p_{2n+1}v=p_nr_1 \label{e8}\\
\Lambda & : & p_n\Lambda\sigma_m=\delta_{n,m} ; \Gamma=\Lambda^{-1}\label{e9} \, (\rm exists\,\,by\,\,hypothesis)\\
 & & r_k(\Lambda\times\Lambda)s_l=\delta_{k,l}\Lambda,\,\, \rm{i.e.}\,\,r_k(\Lambda\times\Lambda)s_l=0\,\,{\rm for}\,\,l\neq k\,\,{\rm and}\,\,r_k(\Lambda\times\Lambda)s_k=\Lambda. \label{e10}
\end{eqnarray}
with $\delta_{k,l}=0$ if $k\neq l$ and $\delta_{k,k}=id$. 
We see that the diagram is commutative: \\
Diagram (1):
\begin{eqnarray*}
p_{2n}v\eta\Pi & = & p_nr_0\eta\Pi \;\;\;({\rm by}\, (\ref{e7})) \\
 & = & p_n\Pi \;\;\; ({\rm by}\, (\ref{e1})) \\
 & = & 1_A \;\;\; ({\rm by}\, (\ref{e4})) \\
p_{2n+1}v\eta\Pi & = & p_nr_1\eta\Pi \;\;\;({\rm by}\, (\ref{e8})) \\
 & = & p_n\Pi \;\;\; ({\rm by}\, (\ref{e1})) \\
 & = & 1_A \;\;\;({\rm by}\, (\ref{e4}))
\end{eqnarray*}
showing that $p_nv\eta\Pi=1_A=p_n\Pi$ for all $n\in \NN$, thus $v\eta\Pi=\Pi$ by the universal property of the product. A similar argument shows that diagram (3) is commutative. For diagram 2, we have, for $k,l\in \{0,1\}$ and $n\in \NN$:
\begin{eqnarray*}
p_{2n+k}v(\Lambda\times\Lambda)u\sigma_{2n+l} & = & p_nr_k(\Lambda\times\Lambda)s_l\sigma_m \;\;\;({\rm by}\,{\rm equations}\,(\ref{e5})-(\ref{e8})) \\
 & = & \delta_{k,l}\cdot p_n\Lambda\sigma_m \;\;\;({\rm by}\,{\rm equation}\,(\ref{e10})) \\
 & = & \delta_{nm}\delta_{kl} \;\;\; ({\rm by}\,{\rm equation}\,\ref{e9}) \\
 & = & \delta_{2n+k,2m+l}
\end{eqnarray*}
thus $p_nv(\Lambda\times\Lambda)u\sigma_m=\delta_{n,m}=p_n\Lambda\sigma_m$ by ($\ref{e9}$) so $\Lambda=v(\Lambda\times\Lambda)u$. Also, we can see that $u,v$ are isomorphisms. Indeed, take $h_0,h_1\in{\rm Hom}(A^{\NN},A^{\NN})$ the unique morphisms such that $p_nh_0=p_{2n}$ and $p_nh_1=p_{2n+1}$ and $v'\in{\rm Hom}(A^{\NN},A^{\NN}\times A^{\NN})$ the only morphism satisfying $r_0v'= h_0$ and $r_1v'=h_1$. Then $v'$ is the only morphism having 
\begin{eqnarray}
p_nr_0v' & = & p_{2n} \label{e11}\\
p_nr_1v' & = & p_{2n+1} \label{e12}
\end{eqnarray}
Then $p_{2n}vv'=p_nr_0v'=p_{2n}$ and $p_{2n+1}vv'=p_nr_1v'=p_{2n+1}$ by $\ref{e7}, \ref{e8}, \ref{e11}, \ref{e12}$ showing that $vv'=1$ and also $p_nr_0v'v=p_{2n}v=p_nr_0$ and $p_nr_1v'v=p_{2n+1}v=p_nr_1$ showing that $v'v=1$. By a similar argument we get that $u$ is an isomorphism, thus $\Lambda=v(\Lambda\times\Lambda)u$ implies that $u\Gamma v=\Gamma\times\Gamma$. Then $\varepsilon(\Gamma\times\Gamma)\eta=\Gamma+\Gamma$ and by the commutative diagrams (1), (2) and (3) we have $\Sigma\Gamma\Pi+\Sigma\Gamma\Pi=\Sigma(\Gamma+\Gamma)\Pi=\Sigma\varepsilon(\Gamma\times\Gamma)\eta\Pi=\Sigma\varepsilon u\Gamma v\eta\Pi=\Sigma\Gamma\Pi$.\\
Now consider the following diagram:
\begin{diagram}
A & \rTo^{\eta'} & A\times A & \rTo^{1_A\times\Pi} & A\times A^{\NN} & \pile{\rTo^{1_A\times\Gamma}\\ \lTo_{1_A\times\Lambda}} & A\times A^{(\NN)} & \rTo^{1_A\times\Sigma} & A\times A & \rTo^{\varepsilon'} & A \\
\dEqual & & (1)' & & \dTo^{g}\uDots_{g'} & (2)' & \uTo^{h}\dDots_{h'} & & (3)' & & \dEqual \\
A & & \rTo_\Pi & & A^{\NN} & \pile{\rTo^\Gamma\\ \lTo_\Lambda} & A^{(\NN)} & & \rTo_\Sigma & & A \\
\end{diagram}
where $\eta'$ and $\varepsilon'$ are the canonical morphisms, and $g$ is defined as follows: if $x_0, x_1$ are the canonical projections of $A\times A^\NN$ then $p_0g=x_0, p_{n+1}g=p_nx_1$ and $g'$ such that $x_0g'=p_0, p_nx_1g'=p_{n+1}$. Then $g$ and $g'$ are inverse isomorphisms and $g$ comes from the canonical equivalence of sets $\NN\sim \NN\cup \{*\}$ which gives an isomorphism $A\times A^{\NN}\simeq A^{\NN}$. $h$ and $h'$ are defined similarly. By computations similar to the ones above we get that diagrams $(1)', (2)', (3)'$ are commutative. Then $g(1_A\times\Lambda)h=\Lambda$ so $1_A\times\Gamma=h\Gamma g$ and then 
\begin{eqnarray*}
 1_A+\Sigma\Gamma\Pi & = & \varepsilon'(1_A\times\Sigma\Gamma\Pi)\eta' \\
 & = & \varepsilon'(1_A\times\Sigma)(1_A\times\Gamma)(1_A\times\Pi)\eta' \\
 & = & \varepsilon'(1_A\times\Sigma)h\Gamma g(1_A\times\Pi)\eta' \;\;\;({\rm by}\,{\rm diagram}\,(2)') \\
 & = & \Sigma\Gamma\Pi \;\;\;({\rm by}\,{\rm commutative}\,{\rm diagrams}\,(1)'\,{\rm and}\,(3)') \\
\end{eqnarray*}
Now take $f:X\rightarrow Y$ an invertible morphism with inverse $(-f):X\rightarrow Y$. Then $f+(-f)=0$ and taking $A=Y$ in the previous diagrams, denoting by $f'=\Sigma\Gamma\Pi f$ we have $\Sigma\Gamma\Pi=\Sigma\Gamma\Pi+\Sigma\Gamma\Pi$ so we get $f'=f'+f'$ by composing to $f$. Also by composing $f+(-f)=0$ with $\Sigma\Gamma\Pi$ we get that $f'$ is invertible, thus $f'=0$. As $1_A+\Sigma\Gamma\Pi=\Sigma\Gamma\Pi$ we get $f+f'=f'$ by composing to $f$ and then $f=0$ (as $f'=0$). Thus the assertion is proved.
\end{proof}

\begin{theorem}
Suppose $\Cc$ is an $\mathbf{AMon}$ category that has at least one nonzero morphism $f:X\rightarrow Y$ that is invertible in the monoid ${\rm Hom}(X,Y)$. Then the functors $\bigoplus\limits_{i\in I}$ and $\prod\limits_{i\in I}$ are isomorphic if and only if $I$ is finite. 
\end{theorem}
\begin{proof}
It follows by Propositions \ref{p1} and \ref{p2}. The converse is true, as like in an $\mathbf{Ab}$ category, the coproduct of a finite family of objects becomes the product as well (biproduct).
\end{proof}

\begin{example}\label{subizo}
Let $\Cc$ be an $\mathbf{AMon}$ category and assume that the diagonal functor $\Delta:\Cc\rightarrow \Cc^I$ has a left adjoint that is also right adjoint, with $I$ an infinite set. Then products and coproducts indexed by any set of cardinal smaller than $I$ exist and finite products and coproducts  are canonically isomorphic \cite{McL}. Define the category $\Dd$ of the objects of the form $H((X_i)_{i\in I})=((X_i)_{i\in I},(\bigoplus\limits_{i\in I}X_i)^{(\NN)},(\prod\limits_{i\in I}X_i)^{\NN})$ with morphisms $((f_i)_{i\in I},(\bigoplus\limits_{i\in I}f_i)^{(\NN)},(\prod\limits_{i\in I}f_i)^{\NN})$. It is easily seen that $\Dd\hookrightarrow\Cc^I\times\Cc\times\Cc$ is a subcategory and that $\Dd$ is isomorphic to $\Cc^{I}$ by the functor $H$. We show that the functors $S=\bigoplus$ and $P=\prod$ from $\Cc^I\times\Cc\times\Cc=\Cc^{I\sqcup\{s,p\}}$ to $\Cc$ are isomorphic on $\Dd$. 
Consider an isomorphism $\theta:\Cc^\NN\simeq \Cc\times\Cc^\NN$ given by a set equivalence $\NN\sim\{*\}\times\NN$. Then we have a commutative diagram (all the coproduct functors as well as the diagonal functors will be denoted simply $\bigoplus$ and $\Delta$):
\begin{diagram}
\Cc^I & \rTo^{\oplus} & \Cc & \rTo^{\Delta} & \Cc^\NN & \pile{\rTo^\oplus \\ \lDots_\Delta} & & & \Cc \\
& & \dDots^{\Delta} & & \dEqual_{\theta} & & & & \dEqual \\
& & I\times\Cc & \rTo_{\Cc\times \Delta} & \Cc\times\Cc^\NN & \pile{\lDots^{I\times\Delta} \\ \rTo_{I\times\oplus}} & \Cc\times\Cc & \pile{\lDots^{\Delta} \\ \rTo_{\oplus}} & \Cc \\
\end{diagram}
The left hand side of the diagram obviously commutes by the definition of $\theta$. The right part of the diagram commutes when the diagonal functors (dotted arrows) are taken under consideration, i.e. $\theta\Delta=(I\times\Delta)\Delta$. But as the coproduct functors are left adjoints to the diagonal functors, we have that $\oplus(I\times\oplus)$ is left adjoint to $(I\times \Delta)\Delta$ and $\oplus\theta^{-1}$ is left adjoint to $\theta\Delta$ and therefore the functors $(I\times\Delta)\Delta$ and $\oplus\theta^{-1}$ are isomorphic and the right part diagram commutes.\\
This shows that the functor $\Cc^I\ni(X_i)_{i\in I}\rightarrow (\bigoplus\limits_{i\in I}X_i)^{(\NN)}\in\Cc^I$ (i.e. $\oplus\Delta\oplus$) is isomorphic to the functor $\Cc^I\ni(X_i)_{i\in I}\rightarrow \bigoplus\limits_{i\in I}X_i\oplus(\bigoplus\limits_{i\in I}X_i)^{(\NN)}\in\Cc^I$ (i.e. $\oplus(I\times\oplus)(I\times\Delta)\Delta\oplus$), thus the functor $\Cc^I\ni(X_i)_{i\in I}\rightarrow(\bigoplus\limits_{i\in I}X_i)^{(\NN)}\oplus(\prod\limits_{i\in I}X_i)^{\NN}\in\Cc$ is isomorphic to $\Cc^I\ni(X_i)_{i\in I}\rightarrow\bigoplus\limits_{i\in I}X_i\oplus(\bigoplus\limits_{i\in I}X_i)^{(\NN)}\oplus(\prod\limits_{i\in I}X_i)^{\NN}\in\Cc$. We obtain that the functor $S$ is isomorphic to the functor $K:=H((X_i)_{i\in I})\longmapsto (\bigoplus\limits_{i\in I}X_i)^{(\NN)}\oplus(\prod\limits_{i\in I}X_i)^{\NN}$ on $\Dd$. Similarly we can show that $R$ is also isomorphic to $K$ on $\Dd$. As $\Cc^I\times\Cc\times\Cc\simeq \Cc^I$ and $R$ and $S$ are isomorphic, by identifying $\Dd$ with a subcategory of $\Cc^I$ we obtain that the functors $\bigoplus\limits_{i\in I}$ and $\prod\limits_{i\in I}$ are isomorphic on a subcategory of $\Cc^I$ (which is not a full subcategory).

\end{example}

\begin{example}
"Complete Lattice Monoids". Let $\Ll$ be the category whose objects are sets $X$ together with an "operation" $\perp:\Pp(X)\rightarrow X$ which fulfills the general associativity $\perp\bigcup(A_i)_{i\in I}=\perp\{\perp A_i\mid i\in I\}$ and an element $0\in X$ such that for $A\subseteq X$, $\perp A\cup\{0\}=\perp A$ and $\perp\emptyset=0$. For example, complete bounded lattices $(L,\leq,\vee,\wedge,0,1)$ are objects of this category when we take $(X,\perp,0)=(L,\vee,0)$ or $(X,\perp,0)=(L,\wedge,1)$. The morphisms of this category are functions $f:(X,\perp_X,0_X)\rightarrow(Y,\perp_Y,0_Y)$ such that $f(\perp_X A)=\perp_Y f(A)$, $\forall A\subseteq X$ and $f(0_X)=0_Y$. For a family of objects $(X_i,\perp_i,0_i)_{i\in I}$ we define an object on the direct product $(\prod\limits_{i\in I}X_i,\prod,0)$ with $\perp A=(\perp_i p_i(A))_{i\in I}$ for $A\subseteq\prod\limits_{i\in I}X_i$ and $0=(0_i)_{i\in I}$, where $p_i$ is the projection of the direct product set $\prod\limits_{i\in I}X_i$. It is easy to see that $\perp$ and $0$ verifies the axioms of $\Ll$. We can also define injections $\sigma_i:X_i\rightarrow \prod\limits_{k\in I}X_k$, $\sigma_i(x_i)=\delta_{i,k}x_i$, that is, the family having $x_i$ on the $i$-th position and $0_k$ elsewhere. The following are easy to check:\\
(i) $\sigma_i$ and $p_i$ are morphisms in $\Ll$;\\
(ii) $\prod\limits_{i\in I}X_i$ is both a left and right adjoint to the diagonal functor from $\Ll$ to $\Ll^I$, with $\sigma_i$ being the canonical injections of the coproduct and $p_i$ the projections of the product.\\
(iii) The operation $\perp$ allows the introduction of a monoid structure of any ${\rm Hom}(X,Y)$, $(f\perp g)=(x\mapsto f(x)\perp_Y g(x))$ with unit $0=(x\mapsto 0)$.\\
Thus we have an example of an {\bf AMon} category where the product and coproduct functors are isomorphic for an infinite set $I$. Note that as for every two objects $X,Y$ in $\Ll$ every element in the monoid ${\rm Hom}(X,Y)$ is idempotent ($f\perp f=f,\,\forall f$), the only invertible one here is $0$.
\end{example}

\section{The module category case}

Let $R$ be a (noncommutative) ring. All tensor products without specification will be considered over $R$. A coring $\Cc$ is a $R-R$-bimodule endowed with a coaction $\Delta:\Cc\longrightarrow \Cc\otimes\Cc$ and a counit $\varepsilon:\Cc\longrightarrow R$ such that both $\Delta$ and $\varepsilon$ are morphisms of $R$-bimodules and if for $c\in\Cc$ we denote by $c_1\otimes c_2=\Delta(c)$ the tensor representation of $\Delta(c)$ (which is actually a sum of tensors) then
$$c_{11}\otimes c_{12}\otimes c_2=c_1\otimes c_{21}\otimes c_{22} $$
$$c_1\varepsilon(c_2)=\varepsilon(c_1)c_2=c$$
for all $c\in \Cc$. A $\Cc$-comodule $M$ is a right $R$-module together with a right $\Cc$ coaction $\rho:M\longrightarrow M\otimes \Cc$ such that if for any $m\in M$ we denote by $m_0\otimes m_1=\rho(m)$ the tensor representation of $\rho(m)$ (the Sweedler-Heynemann notation) then 
$$m_{00}\otimes m_{01}\otimes m_1=m_0\otimes m_{11}\otimes m_{12}$$
$$m_0\varepsilon(m_1)=m$$
for all $m\in M$. A morphism between two right comodules is an $R$-module morphism $f:M\longrightarrow M'$, such that for all $m\in M$, $f(m)_0\otimes f(m)_1=f(m_0)\otimes m_1\in M'\otimes \Cc$. In this way the category  $\Mm^{\Cc}$ of right comodules over $\Cc$ is introduced. \\
Denote by ${}^{*}\Cc={}_{R}\Hom(C,R)$; then there is a multiplication on ${}^{*}\Cc$, defined by $(f\#g)(c)=g(c_1f(c_2))$, $\forall \,f,g\in {}^{*}C$ and $c\in C$. Then the category of right comodules over $\Cc$ is a subcategory of $\Mm_{{}^{*}\Cc}$: if $M\in \Mm^{\CC}$, the right ${}^{*}\Cc$ module structure on $M$ is given by $m\cdot f=m_0f(m_1)$. \\
Let $I$ be a set, $\Mm_R$ the category of right $R$-modules and $\Mm_R^I$ the direct product category of $\Mm_R$ by itself, indexed by the set $I$. Denote by $\delta$ the diagonal functor from $\Mm_R$ to $\Mm_R^I$, $\delta(M)=(M)_{i\in I}$, a family of copies of $M$ indexed by $I$. 

Let $\Cc$ be a coring. Denote by $F$ the forgetful functor from $\Mm^{\Cc}$ to $\Mm_R$, associating to any $\Cc$-comodule $M$ the right $R$-module $M$ and to a morphism of $\Cc$ comodules $f:M\longrightarrow N$ the morphism of right $R$-modules $f$. For any right $R$-module $N$ we have a right $\Cc$-comodule structure on $N\otimes \Cc$ given by the coaction $n\otimes c\longmapsto n\otimes c_1\otimes c_2$, and for any morphism of right $R$-modules $f:N\longrightarrow N'$ we have a morphism of right $\Cc$-comodules $f\otimes \Cc:N\otimes \Cc\longrightarrow N'\otimes \Cc$. Thus we have a functor $G(-)=(-)\otimes \Cc$ from $\Mm_R$ to $\Mm^{\Cc}$ which is right adjoint to the forgetful functor $F$ (see \cite{BW}). 

We recall a few facts on Frobenius Corings, following \cite{BW}, $\xi27$. 

\begin{definition}
An $R$ coring is called Frobenius iff there is an $R$-bimodule map $\eta:R\longrightarrow C$ and a $\Cc$-bicomodule map $\pi:\Cc\otimes \Cc\longrightarrow \Cc$ such that $\pi(c\otimes \eta(r))=\pi(\eta(r)\otimes c)=c$, $\forall\, c\in \Cc$, equivalently there is an element $e\in \Cc^R=\{c\in \Cc\mid rc=cr,\,\forall\,r\in R\}$ and the bicomodule map $\pi$ such that $\pi(e\otimes c)=\pi(c\otimes e)=c$. $(\pi,e)$ is called a Frobenius system.
\end{definition}

The following theorem (see \cite{BW}) gives the connection between Frobenius functors and Frobenius corings.
\begin{theorem}
For a coring $\Cc$ the following statements are equivalent:
\begin{itemize}
\item[(i)] $\Cc$ is a Frobenius coring.
\item[(ii)] The forgetful functor $F$ is a Frobenius functor.
\item[(iii)] The forgetful functor ${}^{\Cc}\Mm\longrightarrow {}_{R}\Mm$ is a Frobenius functor.
\end{itemize}
\end{theorem}

The following finiteness theorem is a key point of the application: (\cite{BW},$27.9$) 

\begin{theorem}\label{FinFrob}
If $\Cc$ is a Frobenius coring, then $\Cc$ is finitely generated and projective as left and also as right $R$-module. 
\end{theorem}


\begin{example}\label{sweedler}
If $R$ is a ring then it becomes an $R$-coring with the comultiplication $\Delta:R\longrightarrow R\otimes R$, $\Delta(r)=1\otimes r=r\otimes 1$ and $\varepsilon(r)=r$, $\forall\,r\in R$. This is the canonical trivial Sweedler coring associated to the ring $R$. Then the category of right $R$-comodules modules coincides to the category of right $R$-modules; a right $R$-module $M$ is a right $R$-comodule by $\rho_M:M\longrightarrow M\otimes R\simeq M$, $\rho_M(m)=m\otimes 1$. 
\end{example}

\begin{example}
Let $(\Cc_i,\Delta_i,\varepsilon_i)_{i\in I}$ be a family of corings. Then the bimodule $\Cc=\bigoplus\limits_{i\in I}\Cc_i$ becomes a coring with comultiplication $\Delta:\Cc\longrightarrow \Cc\otimes \Cc$, $\Delta(\sum\limits_{i\in I}c_i)=\sum\limits_{i\in I}\Delta_i(c_i)$ and counit $\varepsilon:\Cc\longrightarrow R$, $\varepsilon(\sum\limits_{i\in I}c_i)=\sum\limits_{i\in I}\varepsilon_i(c_i)$. 
\end{example}

\begin{proposition}\label{catizo}
If $(\Cc_i,\Delta_i,\varepsilon_i)_{i\in I}$ is a family of corings and $\Cc=\bigoplus\limits_{i\in I}\Cc_i$ is the coring of the previous example, then the product of categories $\prod\limits_{i\in I}\Mm^{\Cc_i}$ is isomorphic to the category $\Mm^{\Cc}$ of comodules over $\Cc$. 
\end{proposition}
\begin{proof}
A proof a little different from the coalgebra case can be done here (see \cite{DNR}): given a family of comodules $(M_{i},\rho_i)_{i\in I}$, where each $M_i$ is a right $\Cc_i$-comodule, consider the right $R$-module $M=\bigoplus\limits_{i\in I}M_{i}$, which becomes a right $\Cc$-comodule by $\rho:\bigoplus\limits_{i\in I}M_{i}\longrightarrow (\bigoplus\limits_{i\in I}M_{i})\otimes \Cc\simeq \bigoplus\limits_{i,j}M_i\otimes C_j$, $\rho(\sum\limits_{i}m_i)=\sum\limits_i\rho_i(m_i)$. It is easy to see that given a family of morphisms $(f_i)_{i\in I}$, $f_{i}:M_{i}\longrightarrow N_{i}$ in $\Mm^{\Cc_i}$, the right $R$-module morphism $f=\bigoplus\limits_{i\in 
I}f_i$ becomes a morphism of right $\Cc$-comodules. In this way we have a functor $S:\prod\limits_{i\in I}\Mm^{\Cc_i}\longrightarrow \Mm^\Cc$, $S((M_{i})_{i\in I})=\bigoplus\limits_{i\in I}M_{i}$. Let now $(M,\rho)$ be a $\Cc$-comodule and denote by $M_i=\rho^{-1}(M\otimes C_{i})$. We identify the right $R$-module $M\otimes C$ with the direct sum $\bigoplus\limits_{i\in I}M\otimes C_{i}$. Denote by $e_i:\Cc\longrightarrow R$ the bimodule morphisms defined by $e_i(\sum\limits_{j\in J}c_j)=\varepsilon(c_i)=\varepsilon_i(c_i)$, for all $\sum\limits_{i}c_i\in \bigoplus\limits_i\Cc_i$. First notice that $M=\bigoplus\limits_{i\in I}M_i$ as right modules. Let $m\in M$; then $\rho(m)=\sum\limits_{i,j}m_{ij}\otimes c_{ij}$, with $c_{ij}\in \Cc_i$, and $m\cdot e_i=\sum\limits_{j}m_{ij}\varepsilon_i(c_{ij})$. We show that $m\cdot e_i\in M_i$ for all $i\in I$. By the comultiplication property we have $\sum\limits_{i,j}m_{ij\,0}\otimes m_{ij\,1}\otimes c_{ij}=\sum\limits_{i,j}m_{ij}\otimes c_{ij\,1}\otimes c_{ij\,2}$. As $\sum\limits_{j}m_{ij\,0}\otimes m_{ij\,1}\otimes c_{ij}\in M\otimes\Cc\otimes\Cc_{i}$, $\sum\limits_{j}m_{ij}\otimes c_{ij\,1}\otimes c_{ij\,2}\in M\otimes\Cc\otimes\Cc_{i}$ (as $c_{ij}\in C_i$) for all $i$ and $M\otimes\Cc\otimes\Cc=\bigoplus\limits_{i\in I}M\otimes\Cc\otimes\Cc_{i}$, we obtain that $\sum\limits_{j}m_{ij\,0}\otimes m_{ij\,1}\otimes c_{ij}=\sum\limits_{j}m_{ij}\otimes c_{ij\,1}\otimes c_{ij\,2}$ for all $i$. Then 
\begin{eqnarray*}
\rho(m\cdot e_i)& = & (\sum\limits_{j}m_{ij}\varepsilon(c_{ij}))_0\otimes(\sum\limits_{j}m_{ij}\varepsilon(c_{ij}))_1\\
& = & \sum\limits_{j}m_{ij\,0}\otimes m_{ij\,1}\varepsilon(c_{ij})\\
& = & \sum\limits_{j}m_{ij}\otimes c_{ij\,1}\varepsilon_i(c_{ij\,2})\\
& = & \sum\limits_{j}m_{ij}\otimes c_{ij}\in M\otimes\Cc_i
\end{eqnarray*}
showing that $m\cdot e_i\in M_i$. We also have that $m=\sum\limits_{i}\sum\limits_{j}m_{ij}\varepsilon(c_{ij})=\sum\limits_{i}m\cdot e_i\in \sum\limits_i M_i$. Now if $m_i\in M_i$ and $\sum\limits_i m_i=0$, then we have that $m_i\cdot e_j=0$ if $i\neq j$ and $m_i\cdot e_i=m_i$ and $0=\sum\limits_i m_i\cdot e_j=m_j$ for all $j$, showing that the sum $\sum\limits_{i}M_i$ is direct. Thus we have a functor $F$ from $\Mm^{\Cc}$ to $\prod\limits_{i}\Mm^{\Cc_i}$ taking a comodule $(M,\rho)$ to the family of (resp. $\Cc_i$) comodules $(M_i=\rho^{-1}(M\otimes \Cc_i),\rho_i=\rho\vert_{M_i}:M_i\longrightarrow M_i\otimes C_i)_{i\in I}$. If $f:M\longrightarrow N$ is a $\Cc$ comodule morphism then this functor takes $f$ into the family of morphisms $(f_i=f\vert_{ M_i})_{i\in I}$, as a simple computation shows that $f(M_i)\subseteq N_i$ and then $f_i$ becomes a morphism of $\Cc_i$ comodules from $M_i=\rho_M^{-1}(M\otimes \Cc_i)$ to $N_i=\rho_N^{-1}(N\otimes\Cc_i)$ (if $m_i\in M_i$, then $\rho_M(m_i)=m_{i\,0}\otimes m_{i\,1}\in M\otimes \Cc_i$ so $\rho_N(f(m_i))=f(m_{i\,0})\otimes m_{i\,1}\in N\otimes\Cc_i$, so $f(M_i)\subseteq N_i$). The previous computations also show that $S\circ F\simeq Id$. A straightforward argument now shows that the functors $F$ and $S$ define inverse equivalences of categories. 
\end{proof}


\begin{theorem}
The functors $\prod\limits_{i\in I}$ and $\bigoplus\limits_{i\in I}$ from $\Mm_R^I$ to $\Mm_R$ (the product and coproduct of modules) are naturally isomorphic if and only if $I$ is a finite set.
\end{theorem}
\begin{proof}
It is easy to see that denoting by $\tilde{R}$ the Sweedler coring from example \ref{sweedler}, we have $\Mm_{R}=\Mm^{\tilde{R}}$ and then by Proposition \ref{catizo} we have $(\Mm_R)^I\simeq \Mm^{\tilde{R}^{(I)}}$, with the described $R$-coring structure of $\tilde{R}^{(I)}$. Consider the following diagram:
\begin{diagram}
(\Mm_R)^I & \simeq & \Mm^{\tilde{R}^{(I)}} \\
\dTo^{\bigoplus} & \ldTo_F & \\
\Mm_R & & \\
\end{diagram}
where $F$ is the forgetful functor. The diagram is commutative, because taking a family of modules $(M_i)_{i\in I}$ to $S((M_i)_{i\in I})$ (where $S$ is the functor from the proof of Proposition \ref{catizo}) and then applying the forgetful functor gives the same as taking the direct sum of modules $\bigoplus\limits_{i\in I}M_i$. The diagram shows that $\bigoplus$ is Frobenius (thus isomorphic to $\Pi$) if and only if $F$ is Frobenius and this implies that $R^{(I)}$ is a finitely generated $R$-module by Theorem \ref{FinFrob}, which is only possible for a finite set $I$.
\end{proof}

\begin{example}
Following example \ref{subizo}, every family of objects in a module category (or any {\bf AMon} category) $(M_{i})_{i\in I}$ can be completed to the family $(M_i)_{i\in I\sqcup\{*\}}$ that has its product and coproduct isomorphic, where $M_{*}=(\bigoplus\limits_{i\in I}M_i)^{(\NN)}\times(\prod\limits_{i\in I}M_i)^{\NN}$. These examples provide a very large class families modules having their product and coproduct isomorphic.
\end{example}

\begin{example}
Let $(X_i)_{i\in I}$ be a family of nonzero (right) comodules over a coalgebra $C$ such that $\Sigma=\bigoplus\limits_{i\in I}X_i$ is a quasifinite comodule. Then $\bigoplus\limits_{i\in I}X_i=\prod^{C}\limits_{i\in I}X_i$, where $\prod^C\limits_{i\in I}$ is the product in the category of comodules.
\end{example}
\begin{proof}
We have that $\prod^C\limits_{i\in I}X_i={\rm Rat}(\prod\limits_{i\in I}X_i)$, where $\prod\limits_{i\in I}$ is the product of modules. Suppose $x=(x_i)_{i\in I}\in P={\rm Rat}(\prod\limits_{i\in I}X_i)$ and the set $J=\{i\in I\mid x_i\neq 0\}$ is infinite. Then $C^*\cdot x$ is a finite dimensional rational comodule, so it has a finite Jordan-Holder composition series. For each $i\in I$ the canonical projection $p_i:C^*\cdot x\rightarrow C^*\cdot x_i$ is an epimorphism, thus $C^*\cdot x_i$ is a rational module (comodule). Consider $\Ss$ a set of representatives for the simple right modules, and $S_i$ a simple subcomodule of $C^*\cdot x_i$ for every $i\in J$. As $\Sigma$ is quasi finite, we have that the set $\Ss'=\{T\in \Ss\mid \exists\,i\in J\,{\rm such}\,{\rm that}\,S_i\simeq T\}$ is infinite and for each $T\in \Ss'$ choose a $k\in J$ such that $T\simeq S_k$ and denote $K$ the set of these $k$. As for every $T\in \Ss'$ there is a $k$ and a monomorphism $T\hookrightarrow C^*\cdot x_k$, it follows then that every Jordan-Holder composition series of $C^*\cdot x$ contains a simple factor isomorphic to $T$. As $C^*\cdot x$ is finite dimensional it follows that the set $\Ff$ of simple left $C^*$ modules appearing as factors in any Jordan-Holder composition series is finite. But $\Ss'\subseteq \Ff$ which is a contradiction to the fact that $\Ss'$ is infinite. Thus $x\in\bigoplus\limits_{i\in I}X_i$, and then ${\rm Rat}(\prod\limits_{i\in I}X_i)\subseteq\bigoplus\limits_{i\in I}X_i$ and the proof is finished as the converse inclusion is obviously true.
\end{proof}

\vspace{2cm}\begin{center}
\sc Acknowledgment
\end{center}
The author wishes to thank his Ph.D. adviser Professor Constantin N\u ast\u asescu for very useful remarks on the subject as well as for his continuous support throughout the past years. He also wishes to thank the referee for useful remarks on the subject that helped improving the presentation.

\vspace*{1cm} 
\begin{flushright}
\begin{minipage}{148mm}\sc\footnotesize

Miodrag Cristian Iovanov\\
University of Bucharest, Faculty of Mathematics, Str.
Academiei 14,
RO-70109, Bucharest, Romania\\
{\it E--mail address}: {\tt
yovanov@walla.com, yovanov@gmail.com}\vspace*{3mm}

\end{minipage}
\end{flushright}
\end{document}